\documentstyle[12pt,amssymb,fleqn]{article}

\def\beq{\begin{equation}}
\def\eeq{\end{equation}}

\begin{document}

\noindent {\em Letters in Mathematical Physics} {\bf 21}:
215--220, 1991.
\bigskip
\bigskip
\bigskip
\bigskip

\noindent
{\Large $q$-Deformed Orthogonal and Pseudo-Orthogonal}
\medskip

\noindent
{\Large Algebras and Their Representations }

\vspace{10mm}

\noindent
{ A. M.~GAVRILIK \ and \  A. U.~KLIMYK }

\noindent {\em Institute for Theoretical Physics, 03143 Kiev,
Ukraine }

\vspace{4mm}
\noindent
(Received: 25 June 1990)

\vspace{9mm}

\noindent {\bf Abstract}. Deformed orthogonal and
pseudo-orthogonal Lie algebras are constructed which differ from
deformations of Lie algebras in terms of Cartan subalgebras and
root vectors and which make it possible to construct
representations by operators acting according to
Gel'fand--Tsetlin-type formulas. Unitary representations of the
$q$-deformed algebras $U_q({\rm so}_{n,1})$ are found.

\vspace{2mm}
\noindent
{\bf AMS subject classifications (1980)}. 16A58, 16A64, 17B10, 81D99.

\vspace{9mm}

\noindent {\bf 1.} In his Letter [1], M.~Jimbo defined a
$q$-deformation $U_q(g)$ of any simple Lie algebra $g$ by means of
its Cartan subalgebra and root elements. M.~Rosso has shown in [2]
that to every integral highest weight there corresponds an
irreducible finite-dimensional representation of $U_q(g)$. For the
$q$-deformed algebra $U_q({\rm sl}(n,{\Bbb C})),$
finite-dimensional irreducible representations were explicitly
constructed by M.~Jimbo [3] through the $q$-analogue of the
Gelfand--Tsetlin formulas. If one constructs the algebra $U_q({\rm
so}(n,{\Bbb C}))$ according to Jimbo's formulas, then (as well as
for the nondeformed case) it is impossible to derive irreducible
finite-dimensional representations in terms of the formulas of the
Gel'fand--Tsetlin type. In explicitly constructing representations
of the Lie algebras of the orthogonal groups, one uses the
generators $I_{k,k-1}=E_{k,k-1}-E_{k-1,k}$ , where $E_{is}$ is the
matrix with the elements $(E_{is})_{jr}=\delta_{ij}\delta_{sr}$.

The purpose of this Letter is to propose another deformation
 $U_q({\rm so}(n,{\Bbb C}))$ of the orthogonal algebras which allows one
to construct the $q$-analogue of the Gel'fand--Tsetlin formulas
for them. With the help of $*$-operations in $U_q({\rm so}(n,{\Bbb
C}))$ , it is possible to introduce a compact deformed algebra
$U_q({\rm so}_n)$ and pseudo-orthogonal deformed algebras
$U_q({\rm so}_{r,s})$,  $r+s=n$. We derive `unitary'
representations (that is, $*$-representations) of the deformed
Lorentz algebras $U_q({\rm so}_{n,1})$ . It turns out that, unlike
the classical Lie algebra ${\rm so}(n,1)$, for the algebras
$U_q({\rm so}_{n,1})$ there also appears a continuous `unitary'
series (strange series) of representations. Under $q\to 1$, this
series disappears (goes to infinity).
\bigskip

\noindent {\bf 2.} By the quantum algebra $U_q({\rm so}(n,{\Bbb
C}))$, $\ n\ge 3$, we shall mean the complex associative algebra
generated by the elements $I_{i,i-1}$, $\ i=2,\ldots ,n$ , which
obey the relations \beq [I_{i,i-1},I_{j,j-1}] =0  \qquad {\rm if}
\quad \mid {i-j}\mid >1,
                                                         \label{(1)}
\eeq \beq I_{i+1,i}^2I_{i,i-1} -
(q^{1/2}+q^{-1/2})I_{i+1,i}I_{i,i-1}I_{i+1,i} +
I_{i,i-1}I_{i+1,i}^2  =  - I_{i,i-1},                  \label{(2)}
\eeq \beq I_{i+1,i} I_{i,i-1}^2 -
(q^{1/2}+q^{-1/2})I_{i,i-1}I_{i+1,i}I_{i,i-1} + I_{i,i-1}^2
I_{i+1,i}  =  - I_{i+1,i} .              \label{(3)} \eeq
In the
limit $q\to 1$, the algebra $U_q({\rm so}(n,{\Bbb C}))$ reduces to
the universal enveloping algebra $U({\rm so}(n,{\Bbb C}))$ of the
classical complex Lie algebra ${\rm so}(n,{\Bbb C})$. Namely,
relations (2) and (3) transform into the relations
\[
[I_{i+1,i},[I_{i+1,i}, I_{i,i-1}]]  =  - I_{i,i-1}\ ,
\]
\[
[I_{i,i-1},[I_{i,i-1}, I_{i+1,i}]]  =  - I_{i+1,i}\ ,
\]
which, together with relations (1), define the Lie algebra ${\rm
so}(n,{\Bbb C})$.

Several ways exist for introducing a $*$-structure (antilinear
antiautomorphism) into $U_q({\rm so}(n,{\Bbb C}))$. The
$*$-structure
$$
I^*_{i,i-1} = - I_{i,i-1}, \ \ \ \ i=2,\ldots ,n\ , \eqno (4a)
$$
determines in $U_q({\rm so}(n,{\bf C}))$ the compact $q$-deformed
algebra $U_q({\rm so}_n)$. The $*$-structure
$$
I^*_{i,i-1} = - I_{i,i-1}, \ \ i\ne r+1; \ \ \ \ \ I^*_{r+1,r} = 
I_{r+1,r} \eqno (4b)
$$
determines in $U_q({\rm so}(n,{\Bbb
C}))$ the noncompact $q$-deformed algebra $U_q({\rm so}_{r,s})$,
$\ r+s=n$.

Let us remark that our definition of the algebra $U_q({\rm
so}(n,{\Bbb C}))$ agrees (when $n=3$) with a new definition of the
quantum algebra $U_q({\rm su}_2)$ introduced by Witten [4]. More
exactly, he defines  $U_q({\rm su}_2)$ as the algebra $U_q({\rm
so}_3)$. Podle{\'s} and Woronowicz [5] have described the Lorentz
quantum group basing on the Lie algebra ${\rm sl}(2,{\Bbb C})$. We
obtain the description of the quantum Lorentz algebra as the
algebra $U_q({\rm so}_{3,1})$ (being a particular case $n=4$ of
our treatment).

The relation of our algebra $U_q({\rm so}(n,{\Bbb C}))$ with the
corresponding algebra introduced by Jimbo [1], and the structure
of a Hopf algebra in it, will be considered in a forthcoming
article.
\bigskip

\noindent {\bf 3.} By a finite-dimensional representation of the
algebra $U_q({\rm so}_n)$, we mean a homomorphism $T$ from this
algebra into the algebra of linear operators acting on
finite-dimensional space such that $T(a^*)=T(a)^*$, $\ a\in
U_q({\rm so}_n)$. To describe a representation of $U_q({\rm
so}_n)$, it is sufficient to give the action formulas for the
operators $T(I_{i,i-1}),\ i=2,3,\ldots ,n$ , which obey relations
(1)--(3) and the relations $T(I_{i,i-1})^* = - T(I_{i,i-1})$.

Finite-dimensional irreducible representations of the algebra
$U_q({\rm so}_n)$ are characterized by highest weights analogous
to the case of the Lie algebra ${\rm so}_n$ [6], that is, by the
numbers $m_{1n},m_{2n},...,m_{kn}$, where $k$ is the integral part
of $n/2$, which satisfy the same conditions as in the classical
case. Namely, all $m_{in}$ are simultaneously integers or
half-integers and
\[
m_{1n}\ge m_{2n}\ge \ldots \ge m_{kn}\ge 0 \hspace{6mm} {\rm for}
\ \ \ n=2k+1\ ,
\]
\[
 m_{1n}\ge m_{2n}\ge \ldots \ge m_{k-1,n}\ge |m_{kn}|
               \hspace{6mm} {\rm for} \ \ \ n=2k\ .
\]
We denote the set of the numbers $m_{1n}, \ldots , m_{kn}$
(highest weight) by a single symbol ${\bf m}$ and the
corresponding representation by $T_{\bf m}$.

In the carrier space of irreducible representation there exists
the basis of the Gel'fand-Tsetlin type labelled by the patterns
$\alpha$ analogous to the case of the Lie algebra $so_n$. The
labels of the patterns $\alpha$ satisfy `betweenness' conditions
[6]. We denote the $k$th row of the pattern $\alpha$ by $m_{1k},
m_{2k},\ldots \ $.

The action of the operators $T_{\bf m}(I_{i,i-1})$ onto
the basis elements $\vert \alpha \rangle \equiv \alpha$ is
written out by means of the $q$-numbers defined as follows.
If $b$ is a complex number, then the $q$-number $[b]$ is given
by the formula
\[
[b] = \frac{ q^{b/2} - q^{-b/2} }
           { q^{1/2} - q^{-1/2} }\ .
\]
Now the action of the operators $T_{\bf m}(I_{i,i-1})$ is given by
the $q$-analogue of the Gel'fand-Tsetlin formulas
$$
 T_{\bf
m}(I_{2p+1,2p}) \alpha = \sum^p_{j=1}A^j_{2p}(\alpha)\alpha^{2p}_j
-\sum^p_{j=1}A^j_{2p}(\bar\alpha^{2p}_j)\bar\alpha^{2p}_j\ , \eqno
(5)
$$ $$
 T_{\bf m}(I_{2p+2,2p+1}) \alpha = \sum^p_{j=1}B^j_{2p}(\alpha)
\alpha^{2p+1}_j -\sum^p_{j=1}B^j_{2p}(\bar{\alpha}^{2p+1}_j)
\bar{\alpha}^{2p+1}_j +
$$ $$
 +{\rm i} \prod_{r=1}^{p+1} [l_{r,2p+2}] \prod_{r=1}^p [l_{r,2p}]
\left(\prod_{r=1}^p [l_{r,2p+1}] [l_{r,2p+1}-1]\right)^{-1}\
 \alpha .
\eqno (6)
$$
where the pattern $\alpha^k_j$ (correspondingly
$\bar{\alpha}^k_j$) is obtained from $\alpha$ by replacing
$m_{jk}$ by $m_{jk}+1$ (correspondingly, by  $m_{jk}-1$);
$l_{j,2p}= m_{j,2p}+p-j,\ l_{j,2p+1}=m_{j,2p+1}+p-j+1$ and
$$
A^j_{2p}(\alpha) =\left( \frac{[l'_j][l'_j+1]}{[2l'_j][2l'_j+2]}
\frac{\prod_{r=1}^p [l_r+l'_j][|l_r-l'_j-1|]
      \prod_{r=1}^{p-1}[l''_r +l'_j][|l''_r -l'_j-1|] }
     {\prod_{r\ne j}[l'_r+l'_j][|l'_r-l'_j|]
         [l'_r+l'_j+1][|l'_r-l'_j-1|]}\right)^{\frac12}  ,
\eqno(7)
$$  $$
 B^j_{2p}(\alpha) = \left( \frac{\prod_{r=1}^{p+1}
[l_r+l'_j][|l_r-l'_j|]
      \prod_{r=1}^{p}[l''_r +l'_j][|l''_r -l'_j|] }
     {\prod_{r\ne j}[l'_r+l'_j][|l'_r-l'_j|]
         [l'_r+l'_j-1] [|l'_r-l'_j-1|]} \cdot
\frac{1}{[l'_j]^2[2l'_j+1][2l'_j-1]} \right)^{\frac12}  . \eqno(8)
$$
(Note that, for brevity, we have denoted $l_{k,2p+1}$ in formula
(7) by $l_k$, $\ l_{j,2p}$ by $l'_j$, and $l_{i,2p-1}$ by $l''_i$;
and likewise, in (8) $l_{i,2p+2}$ by $l_i$, $\ l_{i,2p+1}$ by
$l'_i$ and $l_{j,2p}$ by $l''_j$.)

The validity of relations (1)--(3) and (4a) for the operators
(5) and (6) is proved directly by calculating matrix elements
for the left and right-hand sides of these relations.
We do not reproduce these calculations here because of their
awkwardness.
\vskip 9pt

\noindent {\bf 4.} By a representation $T$ of the algebra
$U_q({\rm so}_{n-1,1})$ we shall mean a homomorphism of this
algebra into the algebra of linear operators in a Hilbert space
$V$ defined on an everywhere dense subspace $D$, such that the
restriction of $T$ onto the subalgebra $U_q({\rm so}_{n-1})$
decomposes into a direct sum of finite-dimensional representations
of this subalgebra, any of which may occur in the decomposition
with a finite multiplicity. Moreover, we suppose that $D$ contains
(as its subspaces) the spaces of irreducible representations of
$U_q({\rm so}_{n-1})$. If conditions (4b) are fulfilled for the
representation operators $T(I_{i,i-1})$, then the representation
is called unitary. Throughout the text, we require that $q$ is not
equal to a root of unity.

To obtain representations of the algebra $U_q({\rm so}_{n-1,1})$,
$\ m_{1n}$ in the heighest weights ${\bf m}=(m_{1n}, m_{2n},\ldots
, m_{kn})$ ($k$ is an integral part of $n/2$) of irreducible
representations of $U_q({\rm so}_{n})$ is replaced by a complex
number $\sigma$. Then the numbers $m_{2,n-1}, m_{3,n-1}, \ldots ,
m_{s,n-1}$ ($s$ is an integral part of $(n-1)/2$) of highest
weights $(m_{1,n-1}, m_{2,n-1},  \ldots , m_{s,n-1})$ of
irreducible representations of the maximal compact subalgebra
$U_q({\rm so}_{n-1})$ must satisfy the same betweenness conditions
as in the case of the representations $T_{\bf m}$ of the algebra
$U_q({\rm so}_{n})$, whereas for $m_{1,n-1}$ we impose the
restriction $m_{2n}\le m_{1,n-1}< \infty$. The spaces of
representations of $U_q({\rm so}_{n-1,1})$ consist of the
$q$-Gel'fand--Tsetlin bases for irreducible representations of
$U_q({\rm so}_{n-1})$. Therefore, the representations of $U_q({\rm
so}_{n-1,1})$ so introduced are infinite-dimensional. The
operators of these representations act according to formulas
(5)--(8) where now $m_{1n}$ is replaced by $\sigma$. It is shown
by direct evaluation that relations (1)--(3) remain valid for such
representation operators. We denote so obtained
infinite-dimensional representations of the algebra $U_q({\rm
so}_{n-1,1})$ by $T(c,{\bf m})$, where ${\bf m}=(m_{2n},
m_{3n},\ldots , m_{pn})$, $c=\sigma +p-1$ if $n=2p$ and $c=\sigma
+p$ if $n=2p+1$. These representations are nothing but the
$q$-analogue of the corresponding representations of the group
${\rm SO}_0(n-1,1)$, see [7, 8].

It is seen from formulas (5)--(8) that the representations $T(c,{\bf m})$
and $T(1-c,{\bf m})$ of the algebra $U_q({\rm so}_{2p,1})$ and the
representations $T(c,{\bf m})$ and  $T(-c,\tilde{\bf m})$ of the
algebra $U_q({\rm so}_{2p+1,1})$ pairwise coincide (here $\tilde{\bf m}$
is obtained from ${\bf m}$ by change $m_{pn}$ into $-m_{pn}$).

Now we shall consider three cases separately: (i) $q=e^{h}$, $h\in
{\Bbb R}$; (ii) $q=e^{{\rm i}h}$, $h\in {\Bbb R}$ and $q$ does not
coincide with a root of unity; (iii) $q=\exp (h_1+{\rm i}h_2)$,
$h_1,h_2\in {\Bbb R}$. When $q=e^{h}$, $h\in {\Bbb R}$, the
function $f(z)=[z]$, $z\in {\Bbb C}$, is seen to have period $4\pi
{\rm i}/h$. In view of that, the representations $T(c,{\bf m})$
and $T(c+4\pi k{\rm i}/h,{\bf m})$, $k\in {\Bbb Z}$ (${\Bbb Z}$ is
the set of integers), coincide, and the representations $T(c,{\bf
m})$ and $T(c+2\pi k{\rm i}/h,{\bf m})$, $k\in {\Bbb Z}$, are
equivalent. The equivalence operator is equal to $-I$, where $I$
is the unit operator. Hence, in this case, we consider only the
representations $T(c,{\bf m})$ with $0\le {\rm Im}\ c<2\pi /h$.
Likewise, at $q=e^{{\rm i}h}$, $h\in {\Bbb R}$, the
representations $T(c,{\bf m})$, $T(c+4\pi k/h,{\bf m})$, $k\in
{\Bbb Z}$, coincide and the representations $T(c,{\bf m})$,
$T(c+2\pi k/h,{\bf m})$, $k\in {\Bbb Z}$, are equivalent.
Therefore, in this case we consider only the representations
$T(c,{\bf m})$, $0\le {\rm Re}\ c<2\pi /h$.
\bigskip

\noindent
THEOREM. {\it The representation $T(c,{\bf m})$ of
$U_q({\rm so}_{2p,1})$ is irreducible if $c$ is not an integer for
integral $m_{2n},m_{3n},\ldots , m_{pn}$, $n=2p+1$, and is not a
half-integer for half-integral $m_{2n},m_{3n},\ldots , m_{pn}$, or
if one of the number $c$, $1-c$ coincides with one of the numbers
$l_j=m_{jn}+p-j+1$, $j=2,3, \ldots ,p$. The representation
$T(c,{\bf m})$ of the algebra $U_q({\rm so}_{2p+1,1})$ is
irreducible either if $c$ is not an integer for integral
$m_{2n},m_{3n},\ldots , m_{p+1,n}$, $n=2p+2$, and is not a
half-integer for half-integral $m_{2n},m_{3n},\ldots , m_{p+1,n}$
or if $c$ coincides with one of the number $l_j=m_{jn}+p-j+1$,
$j=2,3,\ldots ,p$, or $|c|<|l_p+1|$.}
\bigskip

The proof of this theorem is similar to that of the corresponding theorem
for the group ${\rm SO}_0(n-1,1)$ (see, for example, [8]).

To select unitary representations out of the representations
contained in the theorem, one has to verify for which of them the
last condition of relation (4b) is fulfilled. Direct verification
shows that for $q=e^{h}$, $h\in {\Bbb R}$, the following
representations
\bigskip

(1) $T(c,{\bf m})$, $c={\rm i}t+\frac 12$, $0<t<2\pi /h$, and also
$T({\frac 12},{\bf m})$ for integral ${\bf m}$ (principal unitary
series);

(2) $T(c,{\bf m})$, ${\frac 12}<c<s+1$, where $s$ is an integer
($0\le s\le p-1$). Moreover, $s=0$ if $l_{p-1}>1$, or at $l_p=1$
the number $s$ is equal to the largest number $r$ for which
$l_{p-r+1}=r$ (supplementary series);

(3) $T(c,{\bf m})$, ${\rm Im}\ c=\pi /h$, ${\rm Re}\ c>\frac 12$
(strange series)
\bigskip

\noindent
are unitary for the algebra $U_q({\rm so}_{2p,1})$ and the representations
\bigskip

(1) $T(c,{\bf m})$, $c={\rm i}t$, $0\le t<2\pi /h$
(principal unitary series);

(2) $T(c,{\bf m})$, $0<c<s$, where $s$ is an integer
($0\le s\le p$) such that
$l_{p-r+2}=r-1$ for $r=1,2,\ldots , s$ (supplementary series);

(3) $T(c,{\bf m})$, ${\rm Im}\ c=\pi /h$, ${\rm Re}\ c>0$
(strange series)
\bigskip

\noindent
are unitary for the algebra $U_q({\rm so}_{2p+1,1})$.
If $q=e^{{\rm i}h}$, $h\in {\Bbb R}$, then unitary representations
of supplementary series do not arise. In this case, the defining
condition ${\rm Im}\ c=\pi /h$ is replaced by ${\rm Re}\ c=\pi /h$.
If $q=\exp (h_1+{\rm i}h_2)$, $h_1,h_2\in {\Bbb R}$, then
the algebra $U_q({\rm so}_{n,1})$ possesses no nontrivial
unitary representations.

Besides the representations of the algebra $U_q({\rm so}_{n,1})$
listed above, at $q=e^{h}$, $h\in {\Bbb R}$, there are other series
of unitary representations (for example, $q$-analogue of the discrete
series). These representations are subrepresentations of reducible
representations $T(c,{\bf m})$. Such series of representations are
completely similar to corresponding series of the classical
groups ${\rm SO}_0(n,1)$ (see, for example, [7]) and are
described in our preprint [9].
\bigskip

\noindent
{\bf References}
\medskip

\noindent
1. Jimbo, M., {\it Lett. Math. Phys.} {\bf 10}, 63 (1985).

\noindent
2. Rosso, M., {\it Comm. Math. Phys.} {\bf 117}, 581 (1988).

\noindent
3. Jimbo, M., {\it Lect. Notes Phys.} {\bf 246}, 334 (1986).

\noindent
4. Witten, E., Preprint IASSNS--HEP--89/32, Princeton, N.J., 1989.

\noindent
5. Podle{\'s}, P. and Woronowicz, S. L., Report No. 20 of the
Mittag--Lefler Institute,

1989.

\noindent
6. Gel'fand, I. M. and Tsetlin, M. L., {\it Dokl. Akad. Nauk SSSR} {\bf 71},
1017 (1950).

\noindent
7. Klimyk, A. U. and Gavrilik, A. M., {\it J. Math. Phys.} {\bf 20},
1624 (1979).

\noindent
8. Klimyk, A. U., {\it Matrix elements and Clebsch--Gordan
coefficients of group represen- 

tations}, Nauk. Dumka, Kiev, 1979 (in
Russian).

\noindent
9. Gavrilik, A. M., Kachurik, I. I., and Klimyk A. U., Preprint ITP--90--26E,
Kiev, 

1990.

\end{document}